\documentclass[11pt,a4paper]{amsart}
\begin{document}

\title{Generators of some Ramanujan formulas}

\author{Jes\'{u}s Guillera}

\address{Zaragoza (Spain)}

\email{jguillera@gmail.com}

\keywords{Ramanujan's series, WZ-method}

\date{}

\subjclass[2000]{Primary 33C20}

\maketitle

\begin{abstract}
In this paper we prove some Ramanujan-type formulas for $1/\pi$ but without using the theory of modular forms. Instead we use the WZ-method created by H. Wilf and D. Zeilberger and find some hypergeometric functions in two variables which are second components of WZ-pairs that can be certified using Zeilberger's EKHAD package. These certificates have an additional property which allows us to get generalized Ramanujan-type series which are routinely proven by computer. We call these second hypergeometric components of the WZ-pairs generators. Finding generators seems a hard task but using a kind of experimental research (explained below), we have succeeded in finding some of them. Unfortunately we have not found yet generators for the most impressive Ramanujan's formulas. We also prove some interesting binomial sums for the constant $1/\pi^2$. Finally we rewrite many of the obtained series using pochhammer symbols and study the rate of convergence.
\end{abstract}

\section*{Some generators}
We consider the functions in TABLE 1, (generators of first order). Using Zeilberger's EKHAD package \cite{petkovsek} we can check that all these hypergeometric functions $G_1(n,k)$ are second components of WZ-pairs \cite{zeilberger} and at the same time get the first components $F_1(n,k)=C_1(n,k)G_1(n,k)$, where $C_1(n,k)$ is the certificate \cite{wilf}. We know that $F_1(n,k)$ and $G_1(n,k)$ are related by
\[ F_1(n+1,k)-F_1(n,k)=G_1(n,k+1)-G_1(n,k). \]
But in each case considered we have $F_1(0,k)=0$ and so
\[
\sum_{n=0}^{\infty}
\left[ G_1(n,k+1-G_1(n,k) \right]=\sum_{n=0}^{\infty} \left[ F_1(n+1,k)-F_1(n,k) \right]=-F_1(0,k)=0,
\]
which allows us to write
\[
\sum_{n=0}^{\infty} G_1(n,0)=\sum_{n=0}^{\infty} G_1(n,1)= \sum_{n=0}^{\infty} G_1(n,2)=\sum_{n=0}^{\infty} G_1(n,3)=\cdots=constant.
\]
Even more, if we make the substitution of $(-1)^k$ by $\cos\pi k$ the conditions of Carlson's Theorem \cite{bailey} hold and applying it we conclude that for every real number $k$ we have
\[ \sum_{n=0}^{\infty} G_1(n,k)={\rm constant}. \]
For the functions $G_1(n,k)$ we are using we can get the value of the constant by taking the limit at $k=1/2$. And some Ramanujan's formulas are obtained by plugging $k=0.$ If we now define the function, (generators of second order)
\[ G_2(n,k)=F_1(n+1,n+k)+G_1(n,n+k), \]
then from Zeilberger's theorem we can easily derive that
\[ \sum_{n=0}^{\infty} G_2(n,k)=\sum_{n=0}^{\infty} G_1(n,k)={\rm constant}, \]
and again some Ramanujan's formulas appear by plugging $k=0.$ We show these formulas in TABLE 2 \cite{ramanujan}.
\par We can also define the generators of higher order and of course we have
\[
\sum_{n=0}^{\infty} G_1(n,k)=\sum_{n=0}^{\infty} G_2(n,k)=\sum_{n=0}^{\infty} G_3(n,k)=\sum_{n=0}^{\infty} G_4(n,k)=\cdots={\rm constant}.
\]

\par \textit{Observation 1}: With the indicated generators not only have we succeeded in proving some Ramanujan's formulas in an easy way, but also obtained many new interesting formulas.
\par \textit{Observation 2}: An equivalent to the first generator was given by D. Zeilberger \cite{ekhad} as an example to show this new way of proving Ramanujan's formulas.
\par
We now show explicitly the most interesting generalized series we have found classified according to the type of binomials they use in case $k=0.$

\section*{First group of formulas}
\begin{align}
\sum_{n=0}^{\infty} {(\cos\pi k) (-1)^n \over 2^{6n} 2^{2k}} {{2n \choose n}^3 {2k \choose k}^2 \over 2^{2k} {n-1/2 \choose k} {n+k \choose n}} (4n+1) &= {2 \over \pi}, \nonumber \\
\sum_{n=0}^{\infty} {(\cos\pi k) \over 2^{8n} 2^{2k}} {{2n \choose n}^2 {2k \choose k} {2n+2k \choose n+k} \over 2^{2k}
{n-1/2 \choose k}} (6n+2k+1) &={4 \over \pi}, \nonumber \\
\sum_{n=0}^{\infty} {(\cos\pi k) (-1)^n \over 2^{9n} 2^{3k}} {{2n \choose n} {n+k \choose n} {2n+2k \choose n+k}^2  \over
2^{2k} {n-1/2 \choose k}} (6n+2k+1) &={2 \sqrt 2 \over \pi} \nonumber
\end{align}
and
\[
\sum_{n=0}^{\infty} {(\cos\pi k) (-1)^n \over 2^{12n} 2^{2k}} {{2n \choose n}^2 {2n+2k \choose n+k}^2  {n+k \choose n} \over 2^{2k} {n-1/2 \choose k} {2n+k \choose n}} {84n^2+56nk+52n+4k^2+12k+5 \over 2n+k+1} ={16 \over \pi}.
\]

\section*{Second group of formulas}
\begin{align}
\sum_{n=0}^{\infty} {(\cos\pi k) (-1)^n \over 2^{10n} 2^{2k}} {{2n \choose n}^2 {2k \choose k}^2 {4n \choose 2n} \over 2^{2k} {2n-1/2 \choose k} {n+k \choose n}} (20n+2k+3) &={8 \over \pi}, \nonumber \\
\sum_{n=0}^{\infty} {(\cos\pi k) 3^k \over 2^{8n} 2^{4k} 3^{2n}} {{2n \choose n} {4n \choose 2n} {2k \choose k} {2n+2k \choose n+k} \over 2^{2k} {2n-1/2 \choose k}} (8n+2k+1) &={2 \sqrt 3 \over \pi} \nonumber
\end{align}
and
\[
\sum_{n=0}^{\infty} {(\cos\pi k) (-1)^n 3^k \over 2^{12n} 3^n 2^{4k}} {{2n \choose n} {2n+2k \choose n+k} {4n+2k \choose 2n+k} {n+k \choose n} \over 2^{2k} {n-1/2 \choose k}} {56n^2+36nk+34n+4k^2+8k+3 \over 2n+k+1} ={16 \sqrt 3 \over 3 \pi}.
\]

\section*{Experimental method for finding generators}
As an example I am going to explain how I found the formula
\[
\sum_{n=0}^{\infty} {(\cos\pi k) 3^k \over 2^{8n} 2^{4k}3^{2n}} {{2n \choose n} {4n \choose 2n} {2k \choose k} {2n+2k \choose n+k} \over 2^{2k} {2n-1/2 \choose k}} (8n+2k+1)={2 \sqrt 3 \over \pi}.
\]
We consider the function
\[
B(n,k)={{2n+2k \choose n+k}^p \over {2n \choose n}^p}{{4n+2k \choose 2n+k}^q \over {4n \choose 2n}^q} {n+k \choose n}^r {1 \over {2n-1/2 \choose k}} {4n \choose 2n} {2n \choose n}^2.
\]
Observe that for $k=0$ we have
\[ B(n,0)= {4n \choose 2n} {2n \choose n}^2. \]
We are looking for a generator of the type
\[ G(n,k)=f(k) {B(n,k) \over 2^{8n} 3^{2n}} (an+bk+c), \]
such that
\[ \sum_{n=0}^{\infty} G(n,k)=\sum_{n=0}^{\infty} G(n,0)= {2 \sqrt 3 \over \pi}. \]
We can solve this problem in an experimental way looking for integer relations between
\[ {\sqrt 3 \over \pi}, \quad \sum_{n=0}^{\infty} B(n,k), \quad \sum_{n=0}^{\infty} B(n,k)n. \]
The structure of program is like this
\begin{align}
& FOR \quad p=-3 \quad TO \quad 3 \nonumber \\
& \quad FOR \quad q=-3 \quad TO \quad 3 \nonumber \\
& \qquad FOR \quad r=-3 \quad TO \quad 3 \nonumber \\
& u(k)=INTEGER \, \, RELATION \, \, \left( \quad {\sqrt 3 \over \pi}, \quad \sum_{n=0}^{\infty} B(n,k), \quad \sum_{n=0}^{\infty} B(n,k)n \right). \nonumber
\end{align}
I have written a code in PARI-GP following the above idea.
\par For $k=0$ we have, of course, the corresponding Ramanujan's formula. So we must investigate $u(1), u(2)$, etc. For these computational purposes we use the following simplifications
\[
{{2n+2k \choose n+k} \over {2n \choose n}}={\prod_{j=1}^{2k} (2n+j) \over \prod_{j=1}^{k} (n+j)^2}, \quad \quad \quad \quad \quad  {{4n+2k \choose 2n+k} \over {4n \choose 2n}}={\prod_{j=1}^{2k} (4n+j) \over \prod_{j=1}^{k} (2n+j)^2},
\]
\[ {n+k \choose k}={1 \over k!} {\prod_{j=1}^{k} (n+j)}, \quad \quad \quad \quad \quad {2n-1/2 \choose k}={1 \over k!} {\prod_{j=1}^{k} (2n+1/2-j)}.
\]
If there are still integer relations then we know the function $B(n,k)$ is valid and we can easily find the values
of $a,b$, and $c.$ So, it only rests to know the function $f(k)$, but using the EKHAD package this task is very simple.

\section*{A new kind of Ramanujan's formulas}
I show some examples of a new kind of ramanujan's type formulas \cite{guillera}
\begin{align}
\sum_{n=0}^{\infty} {(\cos\pi k) (-1)^n \over 2^{12n} 2^{4k}} {{2n \choose n}^4 {2n+2k \choose n+k} {2k \choose k}^2 \over 2^{2k} {n-1/2 \choose k} {n+k \choose n}} (20n^2+12kn+8n+2k+1) &={8 \over \pi^2}, \nonumber \\
\sum_{n=0}^{\infty} {\cos\pi k \over 2^{16n} 2^{4k}} {{2n \choose n}^4 {4n \choose 2n} {2k \choose k}^3 \over 2^{2k} {2n-1/2 \choose k}[{n+k \choose n}^2} (120n^2+84kn+34n+10k+3) &={32 \over \pi^2} \nonumber
\end{align}
and
\[
\sum_{n=0}^{\infty} {(\cos\pi k) (-1)^n \over 2^{20n} 2^{4k}} {{2n \choose n}^4 {2n+2k \choose n+k}^3 {n+k \choose n} \over 2^{2k} {n-1/2 \choose k}{2n+k \choose n}^2} {P(n,k) \over (2n+k+1)^2}={128 \over \pi^2},
\]
where
\begin{multline}
P(n,k)=3280n^4+4000n^3+1592n^2+232n+13+4592kn^3+3816kn^2+ \\ 884nk+62k+ +1008nk^2+336nk^3+2160n^2k^2+92k^2+40k^3. \nonumber
\end{multline}

For $k=0$ we have

\begin{align}
\sum_{n=0}^{\infty} {(-1)^n {2n \choose n}^5 \over 2^{12n}} (20n^2+8n+1) &={8 \over \pi^2}, \nonumber \\
\sum_{n=0}^{\infty} {{2n \choose n}^4 {4n \choose 2n} \over 2^{16n}} (120n^2+34n+3) &={32 \over \pi^2}, \nonumber \\
\sum_{n=0}^{\infty} {(-1)^n {2n \choose n}^5 \over 2^{20n}} (820n^2+180n+13) &={128 \over \pi^2}. \nonumber
\end{align}
The last one adds roughly three digits per term.

\section*{Some more series}
We sum the following series either by comparing them to other series in this paper at the value $k=0$, or directly by taking the limit as $k \to \infty$. It would also be interesting to evaluate them at $k=1/2$.
\begin{align}
\sum_{n=0}^{\infty} {1 \over 2^{8n} 2^{4k}} {{2n \choose n}^3 {2k \choose k}^2 \over {n+k \choose n}^2} (6n+4k+1) &={4 \over \pi}, \nonumber \\
\sum_{n=0}^{\infty} {(-1)^n \over 2^{6n} 2^{4k}} {{2n \choose n}^2 {2n+2k \choose n+k} {2k \choose k} \over {n+k \choose n}} (4n+2k+1) &={2 \over \pi}, \nonumber \\
\sum_{n=0}^{\infty} {(-1)^n \over 2^{12n} 2^{8k}} {{2n \choose n}^5 {2k \choose k}^4 \over {n+k \choose n}^4} (20n^2+8n+1+24kn+8k^2+4k) &={8 \over \pi^2}. \nonumber
\end{align}

\section*{Rate of convergence}
To study the rate of convergence of the series in this paper and also for other purposes, like getting neat sums when $k$ is a rational number, it is convenient to express the sums in an equivalent form using pochhammer symbols. By doing it to many of the obtained series we have formulas with the following aspect
\begin{align}
\sum_{n=0}^{\infty} {1 \over 2^{2n}} {\left( {1 \over 2} \right)_n \left( {1 \over 2}-k \right)_n \left( {1 \over 2}+k
\right)_n \over n!^2 (1+k)_n} (6n+2k+1) &={4 \over \pi} {2^{2k} \over {2k \choose k}}, \nonumber \\
\sum_{n=0}^{\infty} {(-1)^n \over 2^{3n}} {\left( {1 \over 2}-k \right)_n \left( {1 \over 2}+k \right)_n \left( {1 \over 2}+k \right)_n \over n!^2 (1+k)_n} (6n+2k+1) &={2 \sqrt 2 \over \pi} {2^{3k} \over {2k \choose k}}, \nonumber \\
\sum_{n=0}^{\infty} {1 \over 2^{6n}} { \left( {1 \over 2} \right)_n \left( {1 \over 2}-k \right)_n \left( {1 \over 2}+k
\right)_n \over n!^2 \left( 1+ {k \over 2} \right)_n} { \left( {1 \over 2}+k \right)_n \over \left( {1 \over 2}+ {k \over 2} \right)_n} R(n,k) &={16 \over \pi} {2^{2k} \over {2k \choose k}}, \nonumber
\end{align}
where
\[ R(n,k)={84n^2+56nk+52n+4k^2+12k+5 \over 2n+k+1}. \]
And
\begin{align}
\sum_{n=0}^{\infty} {(-1)^n \over 2^{2n}} {\left( {1 \over 2} \right)_n \left( {1 \over 4}-{k \over 2} \right)_n \left( {3 \over 4}- {k \over 2} \right)_n \over n!^2 (1+k)_n} (20n+2k+3) &={8 \over \pi} {2^{2k} \over {2k \choose k}}, \nonumber \\
\sum_{n=0}^{\infty} {1 \over 3^{2n}} {\left( {1 \over 2}+k \right)_n \left( {1 \over 4}-{k \over 2} \right)_n \left( {3 \over 4}- {k \over 2} \right)_n \over n!^2 (1+k)_n} (8n+2k+1) &={2 \sqrt 3 \over \pi} {2^{4k} \over 3^k {2k \choose k}}, \\ \nonumber
\sum_{n=0}^{\infty} {(-1)^n \over 2^{4n} 3^n} {\left( {1 \over 2}-k \right)_n \left( {1 \over 4}+ {k \over 2}  \right)_n \left( {3 \over 4}+{k \over 2}  \right)_n \over  n!^2 \left( 1+ {k \over 2} \right)_n} {\left( {1 \over 2}+k \right)_n \over \left( {1 \over 2}+ {k \over 2} \right)_n} R(n,k) &={16 \sqrt 3 \over 3 \pi} {2^{4k} \over 3^k {2k \choose k}}, \nonumber
\end{align}
where
\[ R(n,k)={56n^2+36nk+34n+4k^2+8k+3 \over 2n+k+1}. \]

\par The following two formulas are sums of hypergeometric series of type $_5F_4$
\begin{align}
\sum_{n=0}^{\infty} {(-1)^n \over 2^{2n}} {\left( {1 \over 2} \right)_n^3 \left( {1 \over 2}-k \right)_n \left( {1 \over 2}+k \right)_n \over n!^3 (1+k)_n^2} (20n^2+12kn+8n+2k+1) &={8 \over \pi^2} {2^{4k} \over {2k \choose k}^2}, \nonumber \\
\sum_{n=0}^{\infty} {1 \over 2^{4n}} {\left( {1 \over 2} \right)_n^3 \left( {1 \over 4}-{k \over 2} \right)_n \left( {3 \over 4}- {k \over 2} \right)_n \over n!^3 (1+k)_n^2} (120n^2+84kn+34n+10k+3) &={32 \over \pi^2} {2^{4k} \over {2k \choose
k}^2}. \nonumber
\end{align}

\par The rate of convergence is roughly the same of the geometric series obtained taking just the first quotients, in the sense that the number of exact digits they give is asymptotically the same.

\par Now I am going to consider one of the above series and show how to get a more accurate behavior of its rate of convergence. The example I choose is
\[
\sum_{n=0}^{\infty} a(k,n)=\sum_{n=0}^{\infty} {1 \over 2^{4n}} {\left( {1 \over 2} \right)_n^3 \left( {1 \over 4}-{k \over 2} \right)_n \left( {3 \over 4}- {k \over 2} \right)_n \over n!^3 (1+k)_n^2} (120n^2+84kn+34n+10k+3).
\]
The number of exact digits when we sum $N$ terms behaves like
\[
\sum_{n=1}^{N} \log_{10} \left[ {a(k,n) \over a(k, n+1)}
\right] \sim {\ln{16} \over \ln {10}} N + \left( {1 \over 2}+3k
\right) {\ln {N} \over \ln {10}}.
\]
Finally I show another formula I have obtained experimentally

\begin{equation} \nonumber
\sum_{n=0}^{\infty} {(-1)^n 3^{3n} \over 2^{9n}} {\left( {1 \over 2}-k \right)_n \left( {1 \over 6}+ {k \over 3}  \right)_n \left( {5 \over 6}+{k \over 3}  \right)_n \over n!^2 \left( 1+ {k \over 2} \right)_n} {\left( {1 \over 2}+k \right)_n \left( {1 \over 2}+ {k \over 3}  \right)_n \over \left( {1 \over 2}+ {k \over 2} \right)_n \left( {1 \over 2} \right)_n} R(n,k)={32 \sqrt 2 \over \pi} {2^{3k} \over {2k \choose k}},
\end{equation}
where
\[ R(n,k)={616n^3+676n^2+214n+15+440n^2k+312nk+46k+88nk^2+36k^2+8k^3 \over (2n+1)(2n+k+1)}. \]

\section*{Table 1}

\begin{center}
\begin{tabular}{|c|c|c|}
  \hline & & \\
  j & $G_1(n,k)$ & CERTIFICATE \\ & & \\
  \hline & & \\ 1 &
  $\displaystyle {(-1)^n (-1)^k \over 2^{6n} 2^{2k}} {{2n \choose n}^3 {2k \choose k}^2 \over 2^{2k} {n-1/2 \choose k} {n+k \choose n}}(4n+1)$
  & $\displaystyle {4n^2 \over (4n+1)(2n-2k-1)}$ \\  &  &  \\
  \hline & & \\ 2 &
  $ \displaystyle {(-1)^k \over 2^{8n} 2^{2k}} {{2n \choose n}^2 {2k \choose k} {2n+2k \choose n+k} \over 2^{2k} {n-1/2 \choose k}}(6n+2k+1)$
  & $\displaystyle {16n^2 \over (6n+2k+1)(2n-2k-1)}$ \\ & & \\
  \hline & & \\ 3 &
  $\displaystyle {(-1)^n (-1)^k \over 2^{9n} 2^{3k}} {{2n \choose n} {n+k \choose n} {2n+2k \choose n+k}^2 \over 2^{2k} {n-1/2 \choose k}}(6n+2k+1)$&
  $\displaystyle {16n^2 \over (6n+2k+1)(2n-2k-1)}$ \\ & & \\
  \hline & & \\ 4 &
  $\displaystyle {(-1)^n (-1)^k \over 2^{10n} 2^{2k}} {{2n \choose n}^2 {2k \choose k}^2 {4n \choose 2n} \over 2^{2k} {2n-1/2 \choose k} {n+k \choose n}}(20n+2k+3)$ &
  $\displaystyle {64n^2 \over (20n+2k+3)(2n-2k-1)}$ \\ & & \\
  \hline & & \\ 5 &
  $\displaystyle {(-1)^k 3^k\over 2^{8n} 2^{4k} 3^{2n}} {{2n \choose n} {4n \choose 2n} {2k \choose k} {2n+2k \choose n+k} \over 2^{2k} {2n-1/2 \choose k}}(8n+2k+1)$&
  $\displaystyle {36n^2 \over (8n+2k+1)(2n-2k-1)}$ \\ & & \\
  \hline
\end{tabular}
\end{center}

\newpage

\section*{Table 2}

\begin{center}
\begin{tabular}{|c|c|c|}
  \hline & & \\
  $j$ & $G_1(n,0) $ & $G_2(n,0)$ \\ & & \\ \hline
  & & \\ $1$ & $\displaystyle \sum_{n=0}^{\infty} {(-1)^n {2n \choose n}^3 \over 2^{6n}}(4n+1)={2 \over \pi} $ &
  $\displaystyle \sum_{n=0}^{\infty} {{2n \choose n}^3 \over 2^{8n}}(6n+1)={4 \over \pi} $ \\ & & \\ \hline
  & & \\ $2$ & $\displaystyle \sum_{n=0}^{\infty} {{2n \choose n}^3 \over 2^{8n}}(6n+1)={4 \over \pi} $ &
  $\displaystyle \sum_{n=0}^{\infty} {(-1)^n {4n \choose 2n} {2n \choose n}^2 \over 2^{10n}}(20n+3)={8 \over \pi} $ \\ & & \\ \hline
  & & \\ $3$ & $\displaystyle \sum_{n=0}^{\infty} {(-1)^n {2n \choose n}^3 \over 2^{9n}}(6n+1)={2 \sqrt 2 \over \pi} $ &
  $\displaystyle \sum_{n=0}^{\infty} { {4n \choose 2n}^2 {2n \choose n} \over 2^{12n}}{48n^2+32n+3 \over 2n+1}={8 \sqrt 2  \over \pi} $ \\ & & \\ \hline
  & & \\ $4$ & $\displaystyle \sum_{n=0}^{\infty} {(-1)^n {4n \choose 2n} {2n \choose n}^2 \over 2^{10n}}(20n+3)={8 \over \pi} $ &
  $\displaystyle \sum_{n=0}^{\infty} {{2n \choose n}^3 \over 2^{12n}}(42n+5)={16 \over \pi} $ \\ & & \\\hline
  & & \\ $5$ & $\displaystyle \sum_{n=0}^{\infty} {{4n \choose 2n}  {2n \choose n}^2 \over 2^{8n} 3^{2n}}(8n+1)={2 \sqrt 3 \over \pi} $ &
  $\displaystyle \sum_{n=0}^{\infty} {(-1)^n {4n \choose 2n} {2n \choose n}^2 \over 2^{12n} 3^n}(28n+3)={16 \sqrt 3 \over 3 \pi} $ \\ & & \\ \hline
\end{tabular}
\end{center}

\end{document}